\documentclass{article}
\usepackage{amsmath}
\usepackage{amsfonts}

\newenvironment{proof}[1][Proof]{\noindent\textbf{#1.} }{\ \rule{0.5em}{0.5em}}
\input{tcilatex}

\begin{document}

\begin{center}
\begin{equation*}
\text{{\large \ On Salem numbers \ which are exceptional units II}}
\end{equation*}

by\ Toufik Za\"{\i}mi

\ 
\end{center}

\textbf{Abstract. }\textit{We show that for any natural number }$n$ \textit{%
satisfying} $n\equiv 4\func{mod}8$ \textit{and} $n\neq 0\func{mod}5,$\textit{%
\ and for any odd integer }$t\geq (n+6)/2$ \textit{there are infinitely many
Salem numbers }$\alpha $\textit{\ of degree }$2t$\textit{\ such that }$%
\alpha ^{n}-1$\textit{\ is a unit. This result, obtained using a
generalization of a construction due to Gross and McMullen }[5],\textit{\
partially completes the main result of \ }[7].

\medskip


\smallskip


\smallskip

\begin{center}
\textbf{1. Introduction}
\end{center}

We continue the investigation of Salem numbers whose powers are exceptional
units. Recall that an algebraic integer $u$ is said to be an exceptional
unit if both $u$ and $u-1$ are units [6], and a Salem number is a real
algebraic integer greater than $1$ whose other conjugates lie inside the
closed unit disc with at least one conjugate lying on the boundary. Some
properties of the set $\mathbb{T}$ of Salem numbers may be found in$\mathbb{%
\ }$[2].

If $S$ \ denotes the minimal polynomial\ of an element $\alpha $ of $\mathbb{%
T},$ then $S$ has two real roots, namely $\alpha _{1}:=\alpha $ and $\alpha
^{-1},$ and the rest, say $\alpha _{2}^{\pm 1},...,\alpha _{t}^{\pm 1},$ lie
on the unit circle. Also, $\deg (S)=2t\geq 4,$ $\alpha $ is a unit, and the
polynomial $T(x):=(x-(\alpha _{1}+\alpha _{1}^{-1}))(x-(\alpha _{2}+%
\overline{\alpha _{2}}))\cdot \cdot \cdot (x-(\alpha _{t}+\overline{\alpha
_{t}}))\in \mathbb{Z}[x],$ called the trace polynomial of $S$ (or of $\alpha
),$ satisfies 
\begin{equation}
S(x)=x^{t}T(x+\frac{1}{x}).  \tag{1}
\end{equation}%
Clearly, $T$ \ is irreducible and has one root greater than $2,$ namely $%
\alpha +\alpha ^{-1}$ (called in [5] a Salem trace number), and $(t-1)\geq 1$
roots belonging to the interval $(-2,2).$ Conversely, an algebraic integer $%
\beta >2$ of degree $t\geq 2$ whose other conjugates lie in $(-2,2)$ is a
Salem trace number associated to a Salem number $\alpha $ of degree $2t,$
via the relation $\beta =\alpha +\alpha ^{-1}.$

In [6], Silverman conducted numerical investigations on the powers of Salem
numbers which are exceptional units, mainly based on the known small
elements of $\mathbb{T}$ from Boyd's lists [3-4].

By considering the problem of the realization of a monic irreducible
polynomial with integer coefficients as the characteristic polynomial of an
automorphism of the even indefinite unimodular lattice, Gross and McMullen
[5] were concerned with Salem numbers $\alpha $ such that $\alpha ^{2}-1$ is
a unit. They called such numbers $\alpha $ unramified Salem numbers, and
they showed in this case that $\deg (\alpha )/2$ must be odd [5, Proposition
3.3]. Moreover, they gave in [5, Theorem 7.3] a construction leading to the
fact that for any odd integer $t\geq 3$ there are infinitely many unramified
Salem numbers of degree $2t.$ The following result, proved in [7], is a
generalization of this construction.

\smallskip \medskip

\bigskip

\textbf{Theorem 1.1 [7] }\textit{Let }$n\geq 1$ \textit{be an odd }(\textit{%
resp. an even}) \textit{natural number, }$t$ \textit{an integer greater than
or equal to} $(n+3)/2$ (\textit{resp.} $t$ \textit{an odd integer greater
than or equal to} $\ (n+4)/2),$ \textit{and } $D$ \textit{a monic polynomial
with integer coefficients and of degree }$t-(n+3)/2$ (\textit{resp.} \textit{%
of degree} $t-(n+4)/2).$ \textit{If the roots of }$D$\textit{\ }(\textit{if
any, i. e., when} $D(x)\neq 1)$ \textit{lie in the interval }$(-2,2),$%
\textit{\ are distinct, and none of them is a root} \textit{of} \textit{the
polynomial} 
\begin{equation*}
C_{n}(x):=\dprod\limits_{j=1}^{(n-1)/2}(x-2\cos (\frac{2j\pi }{n}))\text{ \
\ \ (\textit{resp.} }C_{n}(x):=\dprod\limits_{j=1}^{(n-2)/2}(x-2\cos (\frac{%
2j\pi }{n}))),
\end{equation*}%
\textit{then, for all sufficiently large} \textit{integers }$a,$\textit{\
the polynomial} 
\begin{equation*}
C_{n}(x)(x-2)D(x)(x-a)-1\text{ \ \ \ \ \ (\textit{resp.} }%
C_{n}(x)(x^{2}-4)D(x)(x-a)-1),
\end{equation*}%
\textit{is a trace polynomial }(\textit{of degree }$t)$ \textit{of a Salem
number }$\alpha $\textit{\ such that }$\alpha ^{n}-1$\textit{\ is a unit.}

\bigskip

With the convention that an empty product is equal to one, Theorem 1.1
implies [5, Theorem 7.3] when $n=2.$ In fact, Theorem 1.1 has been used in
[7] to answer the following question, posed in [7] too:

\textit{For which integers }$n\geq 1$\textit{\ and }$t\geq 2$\textit{\ are
there infinitely many Salem numbers }$\alpha $\textit{\ of degree }$2t$%
\textit{\ such that }$\alpha ^{n}-1$\textit{\ is a unit?}

\bigskip

\textbf{Theorem 1.2 [7] }\textit{Let \ }$(n,t)\in \mathbb{N}^{2}.$ \textit{%
Then, there exist infinitely many Salem numbers }$\alpha $\textit{\ of
degree }$2t$\textit{\ such that }$\alpha ^{n}-1$\textit{\ is a unit,
whenever one of the following conditions holds:}

\textit{(i) }$n$\textit{\ is odd and }$t\geq (n+3)/2;$

\textit{(ii) }$n\equiv 2\func{mod}4,$\textit{\ }$t$\textit{\ is odd and }$%
t\geq (n+4)/2;$

\textit{(iii) }$n=2^{s}$\textit{\ for some integer }$s\geq 2,$\textit{\ }$t$%
\textit{\ is odd and }$t\geq (n+6)/2;$

\textit{(iv) }$n\equiv 4\func{mod}8,$ $n\neq 0\func{mod}3,$\textit{\ }$t$%
\textit{\ is odd and }$t\geq (n+6)/2.$

\bigskip

The first value of $n$ for which Theorem 1.2 does not answer the question
above is $n=12.$ The aim of the present note is to prove the following
result.

\bigskip

\textbf{Theorem 1.3 }\textit{Suppose }$n\equiv 4\func{mod}8$ \textit{and} $%
n\neq 0\func{mod}5.$ \textit{Then, for any odd integer }$t\geq (n+6)/2$ 
\textit{there exist infinitely many Salem numbers }$\alpha $\textit{\ of
degree }$2t$\textit{\ such that }$\alpha ^{n}-1$\textit{\ is a unit.}

\bigskip

Clearly, Theorem 1.3 completes partially Theorem 1.2(iv) for a class of
natural numbers $n$ including the value $12.$ As in [7], to prove Theorem
1.3 we use the polynomials $C_{n},$ defined in Theorem 1.1, the first kind
Chebyshev polynomials $t_{k},$ given by relation: 
\begin{equation*}
t_{k}(2\cos \theta )=2\cos k\theta ,\text{ \ \ }\forall \text{ }(k,\theta
)\in \mathbb{N}\times \lbrack 0,\pi ],
\end{equation*}%
and a variant of Theorem 1.1 with some quadratic factors instead the linear
factor $(x-a).$

Some properties of the polynomials $C_{n}$ and $t_{k}$ are collected in the
following section. The proof of Theorem 1.3, postponed to the last section,
is based on some lemmas presented in Section 3. Throughout, when we speak
about conjugates, the norm, the degree and the minimal polynomial of an
algebraic number (resp. about the degree and the reducibility of a
polynomial) without mentioning the basic field, this is meant over $\mathbb{Q%
}.$

Also, an integer means a rational integer, and when we say that "the" $\gcd $
of two polynomials with integers coefficients is equal to one, this means
that these polynomials have no roots in common. All computations are done
using the system Pari [1].

\ 

\begin{center}
\textbf{2. Some properties of the polynomials }$C_{n}$ \textbf{and} $t_{k}$
\end{center}

From its definition the polynomial $C_{n}$ may be regarded as the trace
polynomial of $U_{n}(x):=(x^{n}-1)/(x-1),$ i. e., $%
U_{n}(x)=x^{(n-1)/2}C_{n}(x+1/x)$ when $n$ is odd (resp. of $%
U_{n}(x):=(x^{n}-1)/(x^{2}-1),$ i. e., $U_{n}(x)=x^{(n-2)/2}C_{n}(x+1/x)$
when $n$ is even$).$ Hence, $C_{1}(x)=C_{2}(x)=1,$ $C_{n}(x)\in \mathbb{Z}%
[x],$ the roots of $C_{n}$ (for $n\geq 3)$ are distinct and belong to the
interval $(-2,2),$ and 
\begin{equation}
\gcd (C_{n},C_{m})=1\Leftrightarrow \gcd (n,m)\in \{1,2\},  \tag{2}
\end{equation}%
as asserted in [7, Lemma 2].

In particular, if $n\equiv 0\func{mod}4,$ then $%
C_{n}(-x-1/x)=U_{n}(-x)/(-x)^{(n-2)/2}=-U_{n}(x)/x^{(n-2)/2}=-C_{n}(x+1/x)$
and so $C_{n}(x)$ is an odd function. Thus, $C_{n}(0)=0,$ and $C_{n}$ has $%
n/4-1$ roots in $(-2,0)$ (and their additive inverses in $(0,2)).$

Also, from their definition, the polynomials $t_{k}$ satisfy the identity: 
\begin{equation}
t_{k+2}(x)=xt_{k+1}(x)-t_{k}(x),\text{ \ }\forall k\in \mathbb{N},  \tag{3}
\end{equation}%
with $t_{1}(x)=x$ and $t_{2}(x)=x^{2}-2,$ and the roots of $t_{k}$ are 
\begin{equation*}
2\cos (\frac{\pi }{2k})>2\cos (\frac{\pi }{2k}+\frac{\pi }{k})>\cdot \cdot
\cdot >2\cos (\frac{\pi }{2k}+\frac{(k-1)\pi }{k}).
\end{equation*}%
Hence, $t_{k}(x)\in \mathbb{Z}[x],$ $\deg (t_{k})=k,$ and $t_{k}$ is monic.
Using (3), a simple induction shows that $t_{2k-1}(0)=0,$ $t_{2k}(0)\neq 0,$ 
$t_{2k-1}(x)$ is an odd function, and $t_{2k}(x)$ is an even function; thus
the roots of any polynomial $t_{k}$ are symmetric with respect to $0.$
Recall also, by [7, Lemma 1], that 
\begin{equation}
\gcd (t_{k},C_{n})=1,\text{ \ }\forall (k,n)\in \mathbb{N}^{2}\text{ with }%
n\neq 0\func{mod}4.  \tag{4}
\end{equation}

\bigskip Also we have the following properties of the polynomials $t_{k}$
and $C_{n}.$

\medskip

\ \ \textbf{Lemma 2.1}\textit{\ Let }$(k,n)\in \mathbb{N}^{2}.$ \textit{%
Then, }

\textit{(i)} $C_{4k}=t_{k}C_{2k}$ \textit{and so} $C_{8k}=t_{2k}C_{4k},$

\textit{(ii) } $\gcd (t_{2k},C_{n})=1,$ \textit{whenever} $n\equiv 4\func{mod%
}8.$

\bigskip

\textbf{Proof. }(i) Since the roots of $t_{k}$ are the numbers $2\cos (\frac{%
\pi (1+2j)}{2k}),$ where $j$ runs through $\{0,1,...,k-1\},$ from the
definition of the polynomials $C_{n},$ given in Theorem 1.1, we easily
obtain the desired equality.

(ii) Assume on the contrary that $\gcd (t_{2k},C_{n})\neq 1$ for some
natural numbers $k$ and $n\equiv 4\func{mod}8.$ Then, there is an integer $%
j\in \{1,...,n/2-1\}$ such that $t_{2k}(2\cos (2j\pi /n))=0.$ Hence, $\cos
(4kj\pi /n)=0,$ $4kj\pi /n=\pi /2+\pi l$ \ for some $l\in \mathbb{Z},$ and
this leads immediately to the contradiction $n\equiv 0\func{mod}8.$%
\endproof%

\smallskip \medskip

\textbf{Lemma 2.2}\textit{\ \ Let }$k\in \mathbb{N},$ \textit{and let} $%
r_{k} $ \textit{denote the number of roots of the polynomial }$t_{k}$\textit{%
\ in the interval} $(0,1).$ \textit{Then, }%
\begin{equation*}
r_{k}=\frac{k-\varepsilon }{6},
\end{equation*}%
\textit{where} $\varepsilon \in \{0,1,\pm 2,3,5\}$ \textit{and} $k\equiv
\varepsilon \func{mod}6.$

\textit{In particular,} $r_{4k}$\textit{\ is even if and only if }$k\equiv 0%
\func{mod}3,$\textit{\ and }$r_{2+4k}$\textit{\ is odd if and only if }$%
k\equiv 1\func{mod}3.$

\bigskip

\textbf{Proof. }Let $2\cos ((1+2j)\pi /2k),$ where $j\in \{0,...,k-1\},$ be
a root of the polynomial $t_{k}.$ Then,%
\begin{equation}
2\cos (\frac{(1+2j)\pi }{2k})\in (0,1)\Leftrightarrow \frac{k}{3}-\frac{1}{2}%
<j<\frac{k-1}{2}.  \tag{5}
\end{equation}%
Writing $k=\varepsilon +6l,$ for some integers $l\geq 0$ and $\varepsilon
\in \{0,1,2,3,-2,5\},$ we easily obtain from (5) that $r_{k}=l,$ and so $%
r_{k}=(k-\varepsilon )/6.$ For example, if $\varepsilon =-2$ and $l\geq 1,$
then (5) gives $2l-7/6<j<3l-1$ and $j\in \{2l-1,...,3l-2\}.$ Hence, there
are $(3l-2)-(2l-1)+1=l$ values for the integer $j,$ and this yields that the
polynomial $t_{k}$ has exactly $l$ roots in the interval $(0,1).$ In the
same way we show the first equality in Lemma 2.2 for the other cases.

Since $2k\equiv 0,$ $\pm 2$ $\func{mod}6,$ it follows from the above that $%
r_{2k}=(k-\eta )/3,$ where $\eta =\varepsilon /2\in \{-1,0,1\}$ and $%
2k\equiv \varepsilon \func{mod}6.$ Therefore, $r_{2k}$ is even if and only
if $k-\eta =6m,$ for some integer $m,$ i. e., $r_{2k}\equiv 0\func{mod}%
2\Leftrightarrow k\equiv 0,\pm 1\func{mod}6,$ and hence 
\begin{equation*}
r_{4k}=r_{2(2k)}\equiv 0\func{mod}2\Leftrightarrow 2k\equiv 0\func{mod}%
6\Leftrightarrow k\equiv 0\func{mod}3,
\end{equation*}%
and%
\begin{equation*}
r_{2+4k}=r_{2(1+2k)}\equiv 0\func{mod}2\Leftrightarrow 1+2k\equiv \pm 1\func{%
mod}6\Leftrightarrow k\equiv 0,-1\func{mod}3.\text{%
\endproof%
}
\end{equation*}

\begin{center}
\bigskip

\newpage

\textbf{3. Some lemmas}
\end{center}

From now on suppose $n\equiv 4\func{mod}8.$ Then, $n\equiv 0\func{mod}4$ and
from the last section we have that $\deg (C_{n})=n/2-1\equiv 1\func{mod}2,$ $%
C_{n}(x)$ is an odd function, $C_{n}(0)=0$ and $C_{n}$ has $n/4$ distinct
roots in $(-2,0].$ Fix also a nonnegative integer $k.$ As mentioned in the
introduction, the proof of Theorem 3.1 is a corollary of the following four
lemmas.

\medskip

\textbf{Lemma 3.1} \textit{For any sufficiently large integer }$a$\textit{\
the polynomial }%
\begin{equation*}
C_{n}(x)(x^{2}-4)t_{4k}(x)(x^{2}-ax+1)-1
\end{equation*}%
\textit{is a trace polynomial }(\textit{of degree }$n/2+4k+3$)\textit{\ of a
Salem number }$\alpha $\textit{\ such that }$\alpha ^{n}-1$\textit{\ is a
unit.}

\bigskip

\textbf{Proof.}\textit{\ }As mentioned in Section 2, $t_{4k}(0)\neq 0$ and $%
t_{4k}(x)$ is an even function having $2k$ roots in $(-2,0).$ Suppose $a\geq
3.$ Then, the polynomial $x^{2}-ax+1$ has two roots $\theta _{a}\in (0,1)$
and $1/\theta _{a}\in (a-1,a),$ and so it is irreducible. It follows, by
Lemma 2.1(ii), that the polynomial%
\begin{equation*}
P(x):=C_{n}(x)(x^{2}-4)t_{4k}(x)(x^{2}-ax+1),
\end{equation*}%
of degree $t:=$ $n/2+4k+3,$ is separable and has $n/4+2k+1$ roots in $%
[-2,0]. $

Let $\eta ^{\ast }$ be the smallest positive root of $%
C_{n}(x)(x^{2}-4)t_{4k}(x).$ Because $\theta _{a}=(a-\sqrt{a^{2}-1})/2$
decreases with $a$ to zero, there is a natural number $a_{0}\geq 3$ such
that $\theta _{a}\leq \theta _{a_{0}}<\eta ^{\ast }$ for all $a\geq a_{0}.$

From now on assume that $a\geq a_{0}$ and the roots of $P,$ say $\eta
_{1},\eta _{2},...,$ $\eta _{t},$\ are labelled so that 
\begin{equation*}
-2=\eta _{1}<\cdot \cdot \cdot <\eta _{n/4+2k+1}=0<\eta _{n/4+2k+2}=\theta
_{a}<\cdot \cdot \cdot <\eta _{t-1}=2<\eta _{t}=1/\theta _{a}.\text{ \ \ }
\end{equation*}%
For each $i\in \{1,...,(t-1)/2\}-\{(n+4)/8+k+1\}$ fix an element $\gamma
_{i} $ of the interval $(\eta _{2i-1},\eta _{2i}).$ Then, 
\begin{equation*}
\gamma _{i}-\eta _{1}>0,...,\ \gamma _{i}-\eta _{2i-1}>0,\ (\gamma _{i}-\eta
_{2i})(\gamma _{i}-\eta _{2i+1})\cdot \cdot \cdot (\gamma _{i}-\eta _{t})>0,
\end{equation*}%
and so $P(\gamma _{i})>0.$ Also for $i=(n+4)/8+k+1,$ fix an element $\gamma
_{i}$ of the interval $(\theta _{a_{0}},\eta _{2i})$ (where $\eta _{2i}=\eta
_{n/4+2k+3}=\eta ^{\ast })$ so that $\eta _{2i-1}<\gamma _{i}<\eta _{2i},$
and $P(\gamma _{i})>0.$

Now, for any $i\in \{1,...,(t-1)/2\}$ set 
\begin{equation*}
a_{i}:=\frac{\gamma _{i}^{2}+1}{\left\vert \gamma _{i}\right\vert }+\frac{1}{%
\left\vert \gamma _{i}C_{n}(\gamma _{i})(\gamma _{i}^{2}-4)t_{4k}(\gamma
_{i})\right\vert }.
\end{equation*}%
Clearly, if $a>a_{i}$ then $a\left\vert \gamma _{i}\right\vert -(\gamma
_{i}^{2}+1)>1/\left\vert C_{n}(\gamma _{i})(\gamma _{i}^{2}-4)t_{4k}(\gamma
_{i})\right\vert $ and hence $P(\gamma _{i})=\left\vert P(\gamma
_{i})\right\vert >1.$ By setting $R_{a}(x)=R(x):=P(x)-1,$ we see that 
\begin{equation*}
R(\gamma _{i})>0>-1=R(\eta _{2i-1})=R(\eta _{2i})
\end{equation*}%
and the polynomial $R$ has two roots, say $\beta _{2i-1}$ and $\beta _{2i},$
such that 
\begin{equation*}
\eta _{2i-1}<\beta _{2i-1}<\gamma _{i}<\beta _{2i}<\eta _{2i}.
\end{equation*}%
Therefore, for each integer $a>A:=\max \{a_{0},a_{1},...,a_{(t-1)/2}\}$ the
polynomial $R$ has $(t-1)$ distinct roots belonging to the interval $(-2,2)$
and the remaining root, say $\beta _{a}=\beta ,$ belongs to the interval $%
(1/\theta ,a)\subset (a-1,a)\subset (2,\infty ),$ as $C_{n}(a)t_{4k}(a)\geq
1,$ $a^{2}-4\geq 5$ and $R(1/\theta )=-1<0<4\leq R(a).$

Consequently, $2\leq \deg (\beta )\leq t$ and $\beta $ is a Salem trace
number of degree at most $t.$ From this point we obtain, similarly as in the
proof of [7, Theorem 3], that there is a constant $A^{\prime }\geq A$ such
that for all $a>A^{\prime }$ the polynomial $R$ is a trace polynomial of a
Salem number $\alpha $ satisfying $\alpha +1/\alpha =\beta .$ Finally, if $%
\zeta $ is a root of $x^{n}-1,$ then $P(\zeta +1/\zeta )=0,$ $R(\zeta
+1/\zeta )=-1,$ and so we get, by (1), $S(\zeta )=-\zeta ^{t},$ where $S$
denotes the minimal polynomial of $\alpha .$ Hence, the absolute value of
the resultant $\dprod\limits_{\zeta ^{n}=1}S_{\alpha }(\zeta )$ of the
polynomials $x^{n}-1$ and $S(x)$ is equal to $1,$ and hence $\alpha ^{n}-1$
is a unit.%
\endproof%

\bigskip

\bigskip

\textbf{Lemma 3.2} \textit{Suppose that the polynomial }$C_{n}t_{2+4k}$ 
\textit{has an even number of roots in the interval }$(0,1).$\textit{\ Then,
for any sufficiently large integer }$a$\textit{\ the polynomial }%
\begin{equation*}
C_{n}(x)(x^{2}-4)t_{2+4k}(x)(x^{2}-ax+a-2)-1
\end{equation*}%
\textit{is a trace polynomial \ }(\textit{of degree }$n/2+4k+5$)\textit{\ of
a Salem number }$\alpha $\textit{\ such that }$\alpha ^{n}-1$\textit{\ is a
unit.}

\bigskip

\textbf{Proof.}\textit{\ }We proceed in the same way as in the proof of
Lemma 3.1. Clearly, $t_{2+4k}(0)\neq 0,$ and $t_{2+4k}(x)$ is an even
function having $2k+1$ roots in $(-2,0).$ Suppose $a\geq 3.$ Then, the
polynomial $x^{2}-ax+a-2$ has two roots $\theta _{a}\in (0,1)$ and $\theta
_{a}^{\prime }\in (a-1,a)\subset (2,\infty ),$ and so it is irreducible. It
follows again by Lemma 2.1(ii) that the polynomial%
\begin{equation*}
P(x):=C_{n}(x)(x^{2}-4)t_{2+4k}(x)(x^{2}-ax+a-2),
\end{equation*}%
of degree $t:=$ $n/2+4k+5,$ is separable and has $n/4+2k+1$ roots in $%
[-2,0]. $

Let $\eta ^{\ast }$ be the greatest root of $C_{n}t_{2+4k}$ in $[0,1).$
Because $\theta _{a}=(a-\sqrt{a^{2}-4a+8})/2$ increases with $a$ to $1,$
there is a natural number $a_{0}\geq 3$ such that $\theta _{a}\geq \theta
_{a_{0}}>\eta ^{\ast }$ for all $a\geq a_{0}.$ From now on suppose that $%
a\geq a_{0}$ and the roots of $P,$ say $\eta _{1},\eta _{2},...,$ $\eta
_{t}, $\ are labelled so that 
\begin{equation*}
-2=\eta _{1}<\cdot \cdot \cdot <\eta _{n/4+2k+2}=0<\cdot \cdot \cdot <\eta
_{n/4+2k+r+3}=\theta _{a}<\cdot \cdot \cdot <\eta _{t-1}=2<\eta _{t}=\theta
_{a}^{\prime },\text{ \ \ }
\end{equation*}%
where $r$ denotes the number of roots of $C_{n}(x)t_{2+4k}(x)$ in $(0,1).$
Since $r$ is even, the integer $n/4+2k+r+3$ is also even.

For each $i\in \{1,...,(t-1)/2\}-\{(n+12)/8+k+r/2\}$ fix an element $\gamma
_{i}$ of the interval $(\eta _{2i-1},\eta _{2i}).$ Then, $\gamma _{i}-\eta
_{1}>0,$ $...,$ $\gamma _{i}-\eta _{2i-1}>0,$ $(\gamma _{i}-\eta
_{2i})(\gamma _{i}-\eta _{2i+1})\cdot \cdot \cdot (\gamma _{i}-\eta _{t})>0,$
and so $P(\gamma _{i})>0.$ For $i=(n+12)/8+k+r/2$ choose an element $\gamma
_{i}$ of the interval $(\eta _{2i-1},\theta _{a_{0}})$ (where $\eta
_{2i-1}=\eta _{n/4+2k+r+2}=\eta ^{\ast })$ so that $\eta _{2i-1}<\gamma
_{i}<\eta _{2i}$ and $P(\gamma _{i})>0.$

Also, for each $i\in \{1,...,(t-1)/2\},$ set 
\begin{equation*}
a_{i}:=\frac{\left\vert \gamma _{i}^{2}-2\right\vert }{\left\vert \gamma
_{i}-1\right\vert }+\frac{1}{\left\vert (\gamma _{i}-1)C_{n}(\gamma
_{i})(\gamma _{i}^{2}-4)t_{2+4k}(\gamma _{i})\right\vert }.
\end{equation*}%
It is clear that $a_{i}$ is well defined, since $\eta _{n/4+2k+r+4}\geq 1$
and $\gamma _{i}\neq 1.$ It follows when $a>a_{i}$ that $a\left\vert \gamma
_{i}-1\right\vert -\left\vert \gamma _{i}^{2}-2\right\vert >1/\left\vert
C_{n}(\gamma _{i})(\gamma _{i}^{2}-4)t_{2+4k}(\gamma _{i})\right\vert $ and
hence $P(\gamma _{i})=\left\vert P(\gamma _{i})\right\vert >1.$ From this
point we argue identically as in the proof of Lemma 3.1, and we obtain that
there is an integer $A\geq \max \{a_{0},a_{1},...,a_{(t-1)/2}\},$ depending
only on $C_{n}t_{2+4k},$ such that for all $a>A$ the polynomial $%
R(x):=P(x)-1 $ is a trace polynomial of a Salem number $\alpha $ with $%
\alpha ^{n}-1$ being a unit.%
\endproof%

\bigskip

\textbf{Lemma 3.3 }\textit{Suppose that the polynomial }$C_{n}t_{2+4k}$ 
\textit{has an odd number of roots in the interval }$(0,1),$ $n\neq 0\func{%
mod}5$ \textit{and } $(r_{2+4k}\equiv 1\func{mod}2$ \textit{or} $%
r_{2+4k}\equiv r_{4k}\equiv 0\func{mod}2).$ \textit{Then, for any
sufficiently large integer }$a$\textit{\ the polynomial }%
\begin{equation*}
C_{n}(x)(x^{2}-4)C_{5}(x)t_{4k}(x)(x^{2}-ax+a-2)-1
\end{equation*}%
\textit{is a trace polynomial }(\textit{of degree }$n/2+4k+5$)\textit{\ of a
Salem number }$\alpha $\textit{\ such that }$\alpha ^{n}-1$\textit{\ is a
unit.}

\bigskip

\textbf{Proof.}\textit{\ }Suppose $a\geq 3,$ and set 
\begin{equation*}
P(x):=C_{n}(x)(x^{2}-4)C_{5}(x)t_{4k}(x)(x^{2}-ax+a-2).
\end{equation*}%
Then, $t:=\deg (P)=$ $n/2+4k+5,$ and the polynomial $x^{2}-ax+a-2$ is
irreducible, as it has two non-integer roots $\theta _{a}\in (0,1)$ and $%
\theta _{a}^{\prime }\in (a-1,a)\subset (2,\infty ).$ Also, Lemma 2.1(ii)
together with the relations (2) and (4) yields $\gcd (t_{4k},C_{5})=1,$ $%
\gcd (t_{4k},C_{n})=1$ and ($\gcd (C_{n},C_{5})=1\Leftrightarrow n\neq 0%
\func{mod}5).$

From now on assume $n\neq 0\func{mod}5.$ Then, the polynomial $P$ is
separable and has $n/4+2k+2$ roots in $[-2,0],$ as $C_{5}(x)=(x+(1+\sqrt{5}%
)/2)(x+(1-\sqrt{5})/2)$ has one negative root.

Let $\eta ^{\ast }$ be the greatest root of $C_{n}t_{4k}C_{5}$ in $(0,1).$
Then, $\eta ^{\ast }\geq (\sqrt{5}-1)/2.$ Because $\theta _{a}=(a-\sqrt{%
a^{2}-4a+8})/2$ increases with $a$ to $1,$ there is a natural number $%
a_{0}\geq 3$ such that $\theta _{a}\geq \theta _{a_{0}}>\eta ^{\ast }$ for
all $a\geq a_{0}$ ($a_{0}$ is certainly greater than $3$ since $\theta
_{3}=(3-\sqrt{5})/2<(\sqrt{5}-1)/2\leq \eta ^{\ast }).$ It follows when $%
a\geq a_{0}$ that the roots of $P,$ say $\eta _{1},\eta _{2},...,$ $\eta
_{t},$\ may be labelled so that 
\begin{equation*}
-2=\eta _{1}<\cdot \cdot \cdot <\eta _{n/4+2k+2}=0<\cdot \cdot \cdot <\eta
_{n/4+2k+r+3}=\theta _{a}<\cdot \cdot \cdot <\eta _{t-1}=2<\eta _{t}=\theta
_{a}^{\prime },\text{ \ \ }
\end{equation*}%
where $r$ denotes the number of roots of $C_{n}t_{4k}C_{5}$ in $(0,1).$

Now, we claim that $r$ is even. Indeed, if $r_{2+4k}\equiv 1\func{mod}2$
then we have from the assumption that the number of roots of $C_{n}$ in $%
(0,1)$ is even, and also Lemma 2.2 gives that $r_{4k}\equiv 1\func{mod}2;$
thus $C_{n}t_{4k}$ has an odd number of roots in $(0,1)$ and so $r$ is even.
In a similar way we get when $r_{2+4k}\equiv r_{4k}\equiv 0\func{mod}2$ that 
$r$ is even, as the hypothesis says that the number of roots of $C_{n}$ in $%
(0,1)$ is odd. It follows that the index of $\theta $ (as a root of $P),$
namely $n/4+2k+r+3,$ is even too. From this point the remaining part of the
proof is identical to the one of Lemma 3.2 by replacing $t_{2+4k}$ by $%
C_{5}t_{4k}.$%
\endproof%

\bigskip

\bigskip

\textbf{Lemma 3.4 }\textit{Suppose that the polynomial }$C_{n}$ \textit{has
an odd number of roots in the interval }$(0,1),$ $n\neq 0\func{mod}5,$ $%
r_{2+4k}\equiv 0\func{mod}2$ \textit{and} $r_{4k}\equiv 1\func{mod}2.$ 
\textit{Then, for any sufficiently large integer }$a$\textit{\ the
polynomial }%
\begin{equation*}
C_{n}(x)(x^{2}-4)C_{5}(-x)t_{4k}(x)(x^{2}-ax+a-2)-1
\end{equation*}%
\textit{is a trace polynomial }(\textit{of degree }$n/2+4k+5$)\textit{\ of a
Salem number }$\alpha $\textit{\ such that }$\alpha ^{n}-1$\textit{\ is a
unit.}

\ \ \medskip

\textbf{Proof.} In the same way as in the last proof set 
\begin{equation*}
P(x):=C_{n}(x)(x^{2}-4)C_{5}(-x)t_{4k}(x)(x^{2}-ax+a-2),
\end{equation*}%
where the integer $a$ is at least $3.$ Then, $t:=\deg (P)=$ $n/2+4k+5,$ and
the polynomial $x^{2}-ax+a-2$ is irreducible, since it has two non-integer
roots $\theta _{a}\in (0,1)$ and $\theta _{a}^{\prime }\in (a-1,a)\subset
(2,\infty ).$ Also, Lemma 2.1(ii) gives that $\gcd (t_{4k}(x),C_{5}(-x))=1,$
as $C_{10}(x)=C_{5}(x)C_{5}(-x).$ Further, we claim that 
\begin{equation*}
\gcd (C_{n}(x),C_{5}(-x))\neq 1\Leftrightarrow n\equiv 0\func{mod}5.
\end{equation*}%
Indeed, because $C_{n}(x)$ is an odd function, the polynomial $C_{5}(x)$ is
a factor of $C_{n}(x)$ if and only if so is the polynomial $C_{5}(-x).$
Further, as the quadratic polynomials $C_{5}(x)$ and $C_{5}(-x)$ are
irreducible, we see that $\gcd (C_{n}(x),C_{5}(-x))\in \{1,C_{5}(-x)\},$ and 
$\gcd (C_{n}(x),C_{5}(x))\in \{1,C_{5}(x)\};$ thus $\gcd
(C_{n}(x),C_{5}(-x))\neq 1\Leftrightarrow \gcd (C_{n}(x),C_{5}(x))\neq 1,$ $%
\gcd (C_{n}(x),C_{5}(-x))\neq 1\Leftrightarrow (C_{10}(x)=C_{5}(x)C_{5}(-x)$
is a factor of $C_{n}(x)),$ and the claim follow from the relation (2).

From now on assume $n\neq 0\func{mod}5.$ Then, the polynomial $P$ is
separable and has $n/4+2k+2$ roots in $[-2,0].$ Let $\eta ^{\ast }$ be the
greatest root of $C_{n}(x)t_{4k}(x)$ in $[0,1).$ Then, $\eta ^{\ast }\geq (%
\sqrt{5}-1)/2.$ Because $\theta _{a}=(a-\sqrt{a^{2}-4a+8})/2$ increases with 
$a$ to $1,$ there is a natural number $a_{0}\geq 3$ such that $\theta
_{a}\geq \theta _{a_{0}}>\eta ^{\ast }$ for all $a\geq a_{0}.$ It follows
when $a\geq a_{0}$ that the roots of $P,$ say $\eta _{1},\eta _{2},...,$ $%
\eta _{t},$ may be labelled so that 
\begin{equation*}
-2=\eta _{1}<\cdot \cdot \cdot <\eta _{n/4+2k+2}=0<\cdot \cdot \cdot <\eta
_{n/4+2k+r+3}=\theta _{a}<\cdot \cdot \cdot <\eta _{t-1}=2<\eta _{t}=\theta
_{a}^{\prime },\text{ \ \ }
\end{equation*}%
where $r$ denotes the number of roots of $C_{n}(x)t_{4k}(x)C_{5}(-x)$ in $%
(0,1)$ (recall that $C_{5}(-x)$ has no root in $(0,1)).$ Because the
assumption says that $r$ is even, the index of $\theta $ (as a root of $P),$
namely $n/4+2k+r+3,$ is even too. From this point the remaining part of the
proof is identical to the one of Lemma 3.3 by replacing $C_{5}(x)$ by $%
C_{5}(-x).$%
\endproof%

\bigskip \newpage

\begin{center}
\textbf{4.} \textbf{Proof of Theorem 1.3}
\end{center}

Suppose $n$ $\equiv 4\func{mod}8$ and let $t$ be an odd integer greater than 
$2+n/2.$ Then, $t=3+n/2+2l$ for some non-negative integer $l,$ and Lemma 3.1
yields immediately the desired result when $l$ is even. Also, if $l=2k+1$
and $k$ is a non-negative integer, then $t=5+n/2+4k$ and Lemma 3.2 implies
Theorem 1.3, whenever the polynomial $C_{n}t_{2+4k}$ has an even number of
roots in the interval $(0,1).$ Finally, if $n\neq 0\func{mod}5$ and $%
C_{n}t_{2+4k}$ has an odd number of roots in the interval $(0,1),$ then the
result follows from Lemma 3.3 (resp. from Lemma 3.4) provided $%
r_{2+4k}\equiv 1\func{mod}2$ or $r_{2+4k}\equiv r_{4k}\equiv 0\func{mod}2$
(resp. provided $r_{2+4k}\equiv 0\func{mod}2$ and $r_{4k}\equiv 1\func{mod}%
2).$%
\endproof%

\bigskip

\bigskip

\bigskip

\begin{center}
\textbf{References}
\end{center}

[1] C. Batut, D. Bernardi, H. Cohen and M. Olivier, \textit{User's Guide to
PARI-GP}, Version 2.5.1 (2012).

[2] M. J. Bertin, A. Decomps-Guilloux, M. Grandet-Hugo, M.
Pathiaux-Delefosse and J. P. Schreiber, \textit{Pisot and Salem numbers},
Birkh\"{a}user Verlag Basel, 1992.

[3] D. W. Boyd, \textit{Small Salem numbers}, Duke Math. J. \textbf{44}
(1977), 315-328.

[4] D. W. Boyd, \textit{Pisot and Salem numbers in intervals of the real line%
}, Math. Comp. \textbf{32} (1978), 1244-1260.

[5] B. H. Gross and C. T. McMullen, \textit{Automorphisms of even unimodular
lattices and unramified Salem numbers}, J. Algebra. \textbf{257} (2002),
265-290.

[6] J. H. Silverman, \textit{Exceptional units and numbers of small Mahler
measure,} Exp. Math. \textbf{4} (1995), 69-83.

[7] T. Za\"{\i}mi, \textit{On Salem numbers which are exceptional units }%
(submitted), arXiv:2309.07723 v2 [math. NT], 27 Sep. 2023, 1-12.

\bigskip

\bigskip

Department of Mathematics and Statistics, College of Science

Imam Mohammad Ibn Saud Islamic University (IMSIU)

P. O. Box 90950

Riyadh 11623 Saudi Arabia

Email: tmzaemi@imamu.edu.sa\textit{\ }

\bigskip

\end{document}